\documentclass{amsart}
\usepackage{graphicx}
\usepackage{amssymb}
\usepackage{epstopdf}
\usepackage{nicefrac}
\usepackage{enumerate}
\usepackage{hyperref}
\usepackage{float}
\usepackage{color}
\usepackage{pst-all}
\usepackage{tabularx}

\newtheorem{thm}{THEOREM}[section]

\newtheorem{cor}[thm]{COROLLARY}

\newtheorem{defn}[thm]{DEFINITION}
\newtheorem{ex}[thm]{EXAMPLE}

\newtheorem{lemma}[thm]{LEMMA}
\newtheorem{prob}[thm]{PROBLEM}
\newtheorem{prop}[thm]{PROPOSITION}

\newtheorem{remark}[thm]{REMARK}

\newcommand{\ds}{\displaystyle}

\newcommand{\Aut}{{\rm Aut}} 
\newcommand{\cC}{{\mathcal C}}
\newcommand{\cD}{{\mathcal D}}
\newcommand{\cG}{{\mathcal G}}
\newcommand{\cGF}{\cG_{\F}} 
\newcommand{\cH}{{\mathcal H}}
\newcommand{\cO}{{\mathcal O}}
\newcommand{\CO}{{\rm CO}} 
\newcommand{\cP}{{\mathcal P}}
\newcommand{\cS}{{\mathcal S}}
\newcommand{\cU}{{\mathcal U}}
\newcommand{\diam}{{\rm diam}} 
\newcommand{\dX}{d_{\fX}} 
\newcommand{\e}{{\varepsilon}} 
\newcommand{\F}{{\mathcal F}}

\newcommand{\Fix}{{\rm Fix}} 
\newcommand{\fS}{{\mathfrak{S}}}
\newcommand{\fX}{{\mathfrak{X}}}
\newcommand{\fY}{{\mathfrak{Y}}}
\newcommand{\Homeo}{{\rm Homeo}} 
\newcommand{\Iso}{{\rm Iso}} 
\newcommand{\mS}{{\mathbb S}}
\newcommand{\mZ}{{\mathbb Z}}
\newcommand{\vp}{{\varphi}}
\newcommand{\whalpha}{\widehat{\alpha}}
\newcommand{\wtalpha}{\widetilde{\alpha}}
\newcommand{\wtphi}{\widetilde{\phi}}
\newcommand{\wtM}{\widetilde{M}}
\newcommand{\G}{\Gamma}

\input xy
\xyoption {all}

\begin{document}

\title{Orbit equivalence and   classification of weak solenoids}
 
\author{Steven Hurder}
 \author{Olga Lukina}
 \email{hurder@uic.edu, olga.lukina@univie.ac.at}
\address{SH: Department of Mathematics, University of Illinois at Chicago, 322 SEO (m/c 249), 851 S. Morgan Street, Chicago, IL 60607-7045}

\address{OL: Faculty of Mathematics, University of Vienna, Oskar-Morgenstern-Platz 1, 1090 Vienna, Austria}

\thanks{Version date: March 6, 2018; rev. August 8, 2019}



  \begin{abstract}
  In this work, we study  minimal equicontinuous actions which are locally quasi-analytic. 
  The first main result shows that for   minimal equicontinuous actions which are locally quasi-analytic,  continuous orbit equivalence of the actions  implies return equivalence. This generalizes results of Cortez and Medynets, and of Li.
  The second main result    is     that if $G$ is a finitely-generated, virtually nilpotent group, then every minimal equicontinuous action    by $G$ is locally quasi-analytic. 
  As an application, we  show that the homeomorphism type of a nil-solenoid is determined by the virtual topological full group of its monodromy action.  
  \end{abstract}

\maketitle

% \tableofcontents

\section{Introduction and main results}\label{sec-intro}

The   works of McCord \cite{McCord1965} and Schori \cite{Schori1966} in the 1960's, and the   work of Fokkink and Oversteegen in \cite{FO2002}, studied the self-homeomorphisms    of weak solenoids. The    authors'  work in \cite{CHL2017} studied the more general problem of the classification, up to homeomorphism,  for weak solenoids, where it  was shown that for   certain classes of weak solenoids,     two weak solenoids are homeomorphic  if and only if  their global monodromy actions are return equivalent.

The global monodromy of a weak solenoid is a minimal equicontinuous Cantor action of a finitely generated group, so the classification problem for weak solenoids motivates the study of invariants for such actions.   The purpose of this work is to consider the topological full group as such an invariant, and to find  conditions on the actions which imply   that they are return equivalent, and hence provide a solution to the classification problem. 
For the class of nil-solenoids, Theorem~\ref{thm-main3} gives  such a solution. 
Before stating our main results precisely, we first recall a few basic concepts.

Let $G$ be a countably generated discrete group, and $\Phi \colon G \to {\Homeo}(\fX)$ be   an action of $G$ on a topological space $\fX$. We also denote the action by $(\fX,G,\Phi)$, and write $g\cdot x$ for $\Phi(g)(x)$. We say that $\Phi$  is a \emph{Cantor action} if $\fX$ is a Cantor space. 

The orbit of a point $x \in \fX$ is the subset $\cO(x) = \{g \cdot x \mid g \in G\}$. 
The action is \emph{minimal} if  for all $x \in \fX$, its   orbit $\cO(x)$ is dense in $\fX$.
The \emph{isotropy group} of $x \in \fX$ is the subgroup 
\begin{equation}\label{eq-isotropyx}
G_x = \{ g \in G \mid g \cdot x = x\} \ . 
\end{equation}

 The action $(\fX,G,\Phi)$ is \emph{effective} if the homomorphism $\Phi$ is assumed to have trivial kernel.
  It  is \emph{free} if for all $x \in \fX$ and $g \in G$,   $g \cdot x = x$ implies that $g = e$, where $e \in G$ denotes the identity of the group. 
Introduce the \emph{isotropy set}
\begin{equation}\label{eq-isotropy}
 \Iso(\Phi) = \{ x \in \fX \mid \exists ~ g \in G ~ , ~ g \ne id ~, ~g \cdot x = x    \} = \bigcup_{e \ne g \in G} \ \Fix(g) \ , 
\end{equation}
where $\Fix(g) = \{x \in \fX \mid g \cdot x = x \}$.
The action $(\fX,G,\Phi)$ is said to be \emph{topologically free}  \cite{BoyleTomiyama1998,Li2015,Renault2008} if the set $\Iso(\Phi) $ is meager in $\fX$. In particular, this means that $\Iso(\Phi)$ has empty interior. 
Note that if $e \ne g \in G$ and $\Phi(g)$ acts trivially on $\fX$, then $ \Iso(\Phi) = \fX$, and thus a topologically free action must be effective.

The action  $(\fX,G,\Phi)$ is \emph{equicontinuous} with respect to a metric $\dX$ on $\fX$, if for all $\e >0$ there exists $\delta > 0$, such that for all $x , y \in \fX$ and $g \in G$ we have
 \begin{equation}\label{eq-equic}
 d_X(x,y) < \delta \qquad \Longrightarrow \qquad d_X(g \cdot x, g \cdot y) < \e .
 \end{equation}
Note that the definition is independent of the choice of the metric $\dX$ on $\fX$.

For a Cantor space $\fX$, let $\CO(\fX)$ denote the collection  of all clopen (compact open) subsets of $\fX$.
Note that for $\phi \in \Homeo(\fX)$ and    $U \in \CO(\fX)$, the image $\phi(U) \in \CO(\fX)$.

We say that $U \subset \fX$  is \emph{adapted} to the action $(\fX,G,\Phi)$ if $U$ is a   non-empty clopen subset, and for any $g \in G$, 
$\Phi(g)(U) \cap U \ne \emptyset$ implies that $\Phi(g)(U) = U$. It follows that
 \begin{equation}\label{eq-adapted}
G_U = \left\{g \in G \mid \vp(g)(U) \cap U \ne \emptyset  \right\}  
\end{equation}
is a subgroup of finite index in $G$, called the \emph{stabilizer} of $U$. 

Denote the restricted action on $U$ by $(U,G_U,\Phi_U)$, and   introduce the \emph{restricted holonomy group}: 
\begin{equation}
\cH_U = image \, \{ \Phi_U \colon G_U \to \Homeo(U)\} \ .
\end{equation}

 \begin{defn}\label{def-re}
For $i=1,2$, let $(\fX_i,G_i,\Phi_i)$ be minimal equicontinuous Cantor actions. We say the actions are \emph{return equivalent} if there exists non-empty clopen subsets $U_i \subset \fX_i$ such that $U_i$ is adapted to the action $\Phi_i$, and there is a homeomorphism $h \colon U_1 \to U_2$ whose induced action $h_* \colon \Homeo(U_1) \to \Homeo(U_2)$  restricts to an isomorphism  $h_* \colon \cH_{U_1} \to \cH_{U_2}$.
 \end{defn}
Note that when $U_i = \fX_i$ and both actions are effective, then this definition reduces to the usual notion of conjugacy of the actions, up to the induced group isomorphism $h_* \colon G_1 \cong \cH_{U_1} \to \cH_{U_2} \cong G_2$. In the terminology of \cite{GPS2017}, this is called an \emph{isomorphism} of the actions. For the general case, though, additional assumptions on the group or the action are required to induce an isomorphism of the actions of the groups $G_1$ and $G_2$ from an isomorphism $h_* \colon   \cH_{U_1} \to \cH_{U_2}$.

 \begin{prob}\label{prob-main1}
Find invariants of minimal equicontinuous Cantor actions which are   sufficient to classify the actions up to return equivalence.  
\end{prob}

In this work, we are concerned with the following conjugacy invariant of a Cantor action $(\fX,G,\Phi)$.
\begin{defn}\label{def-TFG}
The topological full group   $[[\fX,G,\Phi]] \subset \Homeo(\fX)$ is the subgroup consisting of homeomorphisms $\phi \colon \fX \to \fX$ such that for all $x \in \fX$, there exists $U \in \CO(\fX)$ with $x \in U$, and $g \in G$ such that their restrictions to $U$ satisfy $\phi |U = \Phi(g) | U$.
\end{defn}
This definition was introduced in the work by Glasner and Weiss \cite[Section~2]{GW1995},   and has been extensively studied for its relation to the classification problem of Cantor actions, and for the algebraic properties of the group itself, as discussed for example in the S\'eminaire Bourbaki survey \cite{deCornulier2014}, and also the works \cite{BezuglyiKwiatkowski2002,BezuglyiMedynets2008,GPS1999,GM2014}. 
The topological full group  is closely related to the notion of continuous orbit equivalence of actions, as will be discussed in Section~\ref{sec-coe}.   This paper   was motivated by the following result of Cortez and Medynets   in  \cite{CortezMedynets2016}:
\begin{thm}\label{thm-CM}
For $i=1,2$, let $G_i$ be a finitely-generated group, and let $(\fX,G_i,\Phi_i)$  be  a free minimal equicontinuous Cantor action on a   Cantor space $\fX$. If $[[\fX,G_1,\Phi_1]] = [[\fX,G_2,\Phi_2]]$ then the actions $\Phi_1$ and $\Phi_2$ are \emph{structurally conjugate}, hence are return equivalent.  
\end{thm}

We give a generalization of Theorem~\ref{thm-CM} for  minimal equicontinuous Cantor actions that are not free, but satisfy a property called `local quasi-analyticity' as given in Definition~\ref{def-LQA}.

\begin{thm}\label{thm-main1}
For $i=1,2$, let $G_i$ be a finitely-generated group, and $(\fX,G_i,\Phi_i)$ be  a  minimal equicontinuous Cantor action which is locally quasi-analytic.  If $[[\fX,G_1,\Phi_1]] = [[\fX,G_2,\Phi_2]]$ then the actions $\Phi_1$ and $\Phi_2$ are return equivalent.
\end{thm}

Our second main result gives a sufficient condition for the group $G$ which implies the LQA property.
 
\begin{thm}\label{thm-main2}
Let $G$ be a Noetherian group.
Then     a minimal equicontinuous action   $(\fX,G,\Phi)$  on a Cantor space $\fX$   is locally quasi-analytic.
\end{thm}

In particular,   a finitely-generated nilpotent group  is Noetherian, so  Theorem~\ref{thm-main1} implies:

 \begin{cor}\label{cor-vpc}
Let  $G$ be a finitely-generated   nilpotent group. Then the topological full group  is a complete invariant up to return equivalence for a minimal equicontinuous   action $(\fX,G,\Phi)$.   If two such actions have spatially isomorphic topological full groups, then the actions are return equivalent.
\end{cor}

This result can be considered as a localized version of   \cite[Theorem~1.3]{Li2015}, but with the Noetherian  hypothesis on $G$  in place of the homological criterion of Li in \cite[Section~5]{Li2015}.

Finally, in Definition~\ref{def-virtualfullgroup}, we introduce  the virtual isomorphism class  $[[[\fX,G,\Phi]]]$ of the topological full group for a minimal equicontinuous action $(\fX,G,\Phi)$, and in Theorem~\ref{thm-VFHrigidity} give an extension of Theorem~\ref{thm-main1} using this invariant of actions.
 As an application, in Section~\ref{sec-nil}  we recall the notion of a nil-solenoid, which is  a weak solenoid whose base  space is a closed nil-manifold and whose global monodromy action is effective.   
Then we prove Theorem~\ref{thm-conjugate}, which implies the following:
\begin{thm}\label{thm-main3}
The   virtual isomorphism class  $[[[\fX,G,\Phi]]]$ for the monodromy of a  nil-solenoid   is a complete invariant of the homeomorphism class of the solenoid. 
  \end{thm}

  Appendix~\ref{sec-examples} contains a collection of examples of minimal Cantor actions which illustrate   properties of   actions as discussed in this work.

\section{Locally quasi-analytic  and Hausdorff actions} \label{sec-ad}

  In this section, we discuss   notions of pointwise and uniform ``regularity'' for topological actions. 
  
  We  first recall the   notion of a locally quasi-analytic action. Haefliger used   in \cite{Haefliger1985} the notion of a quasi-analytic action  on a topological space, in his study of pseudogroups of local isometries on locally connected spaces. 
 The works  \cite{ALC2009,ALM2016} by \'Alvarez-L\'opez, Candel and Moreira-Galicia reformulated   Haefliger's definition for the case of 
 topological actions on  Cantor spaces  as follows:

 \begin{defn} \cite[Definition~9.4]{ALC2009} \label{def-LQA} A topological action       $(\fX,G,\Phi)$  is   \emph{locally quasi-analytic}, or simply   \emph{LQA}, if there exists $\e > 0$ such that for any non-empty open set $U \subset \fX$ with $\diam (U) < \e$,  and  for any non-empty open subset $V \subset U $, and elements $g_1 , g_2 \in G$
 \begin{equation}\label{eq-lqa}
  \text{if the restrictions} ~~ \Phi(g_1)|V = \Phi(g_2)|V, ~ \text{ then}~~ \Phi(g_1)|U = \Phi(g_2)|U. 
\end{equation}
 The action is said to be \emph{quasi-analytic} if \eqref{eq-lqa} holds for $U=\fX$.
\end{defn}

Examples of equicontinuous Cantor actions   $\fX$ which are locally quasi-analytic, but not quasi-analytic,  are given in \cite{DHL2016c,HL2017a}, and  also in Section~\ref{sec-examples}.

The idea of the proof for the following result appeared  in the work   \cite{EMT1977} by Epstein, Millet and Tischler,  in the context of   pseudogroup actions. It  has also appeared in the literature in various alternative formulations, for example as  Proposition~3.6  in \cite{Renault2008} and  Lemma~2.2 in \cite{Li2015}.
\begin{prop}\label{prop-BCT}
An effective Cantor action $(\fX,G,\Phi)$ is quasi-analytic  if and only if it is topologically free.
\end{prop}
\proof
Suppose that the action $\Phi$ is   topologically free, then the isotropy set   $\Iso(\Phi) $ is meager in $\fX$. In particular, this means that $\Iso(\Phi)$ has empty interior. Let $V \subset \fX$ be a non-empty  open set, and suppose that    $g_1,g_2 \in G$ satisfy 
$\Phi(g_1)|V = \Phi(g_2)|V$. Then $\Phi(g_2^{-1}g_1)|V = id|V$ so we must have that $g_2^{-1}g_1 = e \in G$. Thus, $g_1 = g_2$ and hence $\Phi(g_1)  = \Phi(g_2)$, as was to be shown.

Conversely, suppose that the action $\Phi$ is  quasi-analytic and effective, then for each $e \ne g \in G$ the complement of $\Fix(\Phi(g))$ is open and dense in $\fX$, so $\Fix(\Phi(g))$ is a closed nowhere dense set. Thus, $\Iso(\Phi)$ is the countable union of closed nowhere dense sets, hence by the Baire Category Theorem, $\Iso(\Phi)$ must be meager.\endproof

\begin{cor}\label{cor-abelian}
Let  $(\fX,G,\Phi)$ be a minimal topological action. If $G$ is a finitely-generated abelian group, and the action  is effective, then it is quasi-analytic, hence is topologically free.
\end{cor}
 \proof
 Suppose that there exists a non-empty open set $V \subset \fX$ and let $g_1 , g_2 \in G$ be such that the restrictions satisfy $\Phi(g_1)|V = \Phi(g_2)|V$. Let $g = g_2^{-1} \, g_1$ then  the restriction $\Phi(g) | V$ is the identity map. 
 Let $y \in \fX$ then there exists $g_y \in G$ such that $g_y \cdot y \in V$, 
 so $y \in V_y \equiv g_y^{-1} \cdot V$. Since $\Phi(g) = \Phi(g_y^{-1} \, g \, g_y)$, 
 the restriction $\Phi(g) | V_y$ is the identity map. 
 Thus, $\Phi(g)$ is the identity on $\fX$, and as the action is effective, we have $g = id$.
 \endproof

 The LQA property for a group action  $(\fX,G,\Phi)$ can be interpreted in terms of the properties of the germinal groupoid $\cG_{\Phi}$ associated to the action. This groupoid  is fundamental for the study of the $C^*$-algebras these actions generate, as discussed for example by Renault in \cite{Renault1980,Renault2008}.   
 Recall that for $g_1, g_2 \in G$, we say that $\Phi(g_1)$ and $\Phi(g_2)$ are \emph{germinally equivalent} at $x \in \fX$ if $\Phi(g_1) (x) = \Phi(g_2)(x)$, and there exists an open neighborhood $x \in U \subset \fX$ such that the restrictions agree, $\Phi(g_1)|U = \Phi(g_2)|U$. We then write $\Phi(g_1) \sim_x \Phi(g_2)$. For $g \in G$, denote the equivalence class of  $\Phi(g)$ at $x$ by $[g]_x$.
 The collection of germs $\cG(\fX, G, \Phi) = \{ [g]_x \mid g \in G ~ , ~ x \in \fX\}$ is given the sheaf topology, and   forms an \emph{\'etale groupoid} modeled on $\fX$. 
 We recall the following result from Winkelnkemper:
 \begin{prop}\cite[Proposition 2.1]{Winkelnkemper1983a}\label{prop-hausdorff}
The germinal groupoid $\cG(\fX, G, \Phi)$ is Hausdorff at $[g]_x$  if and only if, for all $[g']_x \in \cG(\fX, G, \Phi)$ with $g \cdot x = g' \cdot x = y$, if there exists a sequence   $\{x_i\} \subset \fX$ which converges to $x$ such that $[g]_{x_i} = [g']_{x_i}$ for all $i$, then $[g]_{x} = [g']_{x}$. 
\end{prop}

  Winkelnkemper showed in \cite[Proposition~2.3]{Winkelnkemper1983a} that for a smooth foliation $\F$ of a connected manifold $M$ for which the associated holonomy pseudogroup $\cGF$ is generated by real analytic maps, then $\cGF$ is a Hausdorff space. For Cantor actions, an analogous result holds for the LQA property. 

\begin{prop}\label{prop-hausdorff2}
If an action $(\fX,G,\Phi)$   is   locally quasi-analytic, then  $\cG(\fX, G, \Phi)$ is    Hausdorff.     
 \end{prop}
 \proof
Assume that $\cG(\fX, G, \Phi)$ is  not  Hausdorff.  Then there exists $g \in G$ and $x \in \fX$ such that   $\cG(\fX, G, \Phi)$ is non-Hausdorff at $[g]_x$. By Proposition~\ref{prop-hausdorff},  there exists $[g']_x \in \cG(\fX, G, \Phi)$ with $g \cdot x = g' \cdot x = y$, and a sequence   $\{x_i\} \subset \fX$ which converges to $x$ such that $[g]_{x_i} = [g']_{x_i}$ for all $i$, but $[g]_x \ne [g']_x$. Let  $g'' = g^{-1} g' \in G$,  then $g'' \cdot x = x$ and $[g'']_{x_i} = [id]_{x_i}$ for all $i$, but $[g'']_{x} \ne [id]_{x}$. In particular, the action of $\Phi(g'')$ is not the identity in any open neighborhood of $x$, but there exists a sequence of open sets $x_i \in U_i \subset \fX$ containing $x$ in their limit for which the restriction $\Phi(g'')| U_i = id |U_i$. Hence, there does not exists $\e >0$ such that  $\Phi(g)|U$ is   quasi-analytic  for all open neighborhood $x \in U$ with $\diam \ (U) < \e$. Thus, the action $(\fX,G,\Phi)$   is not  locally quasi-analytic.
\endproof
 
 \begin{remark}
{\rm 
 Suppose that   $(\fX, G, \Phi)$ is a Cantor action such that  $\cG(\fX, G, \Phi)$ is    Hausdorff. 
 Then for each   $x \in \fX$ and $g,g' \in G$ with $g \cdot x = g' \cdot x$ and $[g]_x \ne [g']_x$, there exists an open neighborhood $x \in U(x,g,g') \subset \fX$ such that the set $\{y \in U(x,g,g') \mid g \cdot y = g' \cdot y\}$ has no interior. If there exists $\e > 0$ such that any   open set $U$ with $x \in U$ and $\diam \ (U) < \e$ then $U(x,g,g') = U$ for all $g,g' \in G$,   this is just saying that  $\cG(\fX, G, \Phi)$  locally quasi-analytic. Thus, the  locally quasi-analytic property can be viewed as a ``uniform Hausdorff property'' for $\cG(\fX, G, \Phi)$. 
} \end{remark}

We say that  $x \in \fX$   is     a \emph{non-Hausdorff point} for the action $(\fX,G,\Phi)$ if there exists $g \in G$ such that the germ $[g]_x$ is not Hausdorff in $\cG(\fX, G, \Phi)$. Examples of minimal equicontinuous Cantor actions with non-Hausdorff points are discussed in Section~\ref{sec-examples}. 
The existence of non-Hausdorff points has implications for the algebraic structure of the reduced $C^*$-algebra $C^*(\fX,G,\Phi)$ associated to the action, as discussed for example in  \cite{BCFS2014,Exel2011,Renault2008}.

\section{Equicontinuous  actions} \label{sec-stable}

In this section, we consider some basic properties of  minimal equicontinuous Cantor actions. The main results of the section are  Theorem~\ref{thm-LQA}, which gives an algebraic criterion for    when a group action is locally quasi-analytic, and Corollary~\ref{cor-nilLQA} which applies this criterion to nilpotent groups.
  
 First, we recall the following folklore result and give a sketch of the proof, as the ideas involved are central to the study of equicontinuous actions and used in the discussions that follow. 
 \begin{prop}\label{prop-CO}
Let $G$ be a finitely-generated group. Then a minimal Cantor action $(\fX,G,\Phi)$ is  equicontinuous  if and only if, for the induced action $\Phi_* \colon G \times \CO(\fX) \to \CO(\fX)$, the $G$-orbit of every $U \in \CO(\fX)$ is finite.
\end{prop}
\proof
Suppose that the action is equicontinuous, and let $U_1 \subset \fX$ be a non-empty proper clopen subset. Then $U_2 = \fX - U_1$ is also a non-empty proper clopen subset, and let $\e = \dX(U_1,U_2) > 0$, where $\dX$ is the choice of some metric on $\fX$. Let $\delta > 0$ be such that \eqref{eq-equic} holds for this choice of $\e$, for all $g \in G$.
The iterates of the partition $\{U_1,U_2\}$ of $\fX$ define a closed partition of $\fX$, 
 $$\cC = \bigcap ~ \left\{  g \cdot U_i \mid i =1,2 , ~ g \in G  \right\} \ ,$$
 which is invariant for  the action $\Phi$ by construction.
Let $K \in \cC$ with $x \in K$, and suppose that $g \cdot x \in U_i$. Let $y \in \fX$ satisfy $\dX(x,y) < \delta$ then by \eqref{eq-equic} and the choice of $\e$ we have that $g \cdot y \in U_i$ also. As this holds for all $g \in G$, the set $K$ is open in $\fX$. Thus, $\cC$ is a clopen partition, and by the compactness of $\fX$ it must be finite. Thus, the action $\Phi$  permutes this finite collection, hence there is a   subgroup $G_{\cC} \subset G$ of finite index which fixes every $K \in \cC$. As $U_1$ is the union of sets in $\cC$, the subgroup $G_{\cC}$ also fixes $U_1$, and thus the $G$-orbit of $U_1$ is finite. 

Conversely, assume that the $G$-orbit of every $U \in \CO(\fX)$ is finite. Fix a basepoint $x \in \fX$, and as $\fX$ is a Cantor space, one can choose a descending chain of clopen sets 
$$\fX = V_0 \supset V_1 \supset V_2 \supset \cdots \supset V_{\ell} \supset \cdots \supset \{x\}$$
 whose intersection is $\{x\}$. For each $\ell > 0$, the intersection of the finite collection $\{g \cdot V_{\ell} \mid g \in G\}$ has a clopen set containing the basepoint $x$, which we label $U_{\ell}$, so that $x \in U_{\ell} \subset V_{\ell}$. Then $U_{\ell}$ is an adapted clopen subset with stabilizer group denoted by $G_{\ell} \subset G$. 
 It then follows as in \cite[Appendix A]{DHL2016a} or \cite[Section~2]{CortezMedynets2016}, that the Cantor action $(\fX,G,\Phi)$ is conjugate to the odometer constructed from the group chain $\{G_{\ell} \mid \ell \geq 0\}$, hence    $(\fX,G,\Phi)$ is a minimal equicontinuous action.
\endproof
 
 The above proof that each $U \in \CO(\fX)$ has finite orbit is essentially the same as what was called the ``coding method'' used to study equicontinuous pseudogroup actions in \cite{ClarkHurder2013}, and discussed for group actions  in \cite[Appendix A]{DHL2016a}.
  
We discuss in more detail another key idea in the above proof.   Let $U \subset \fX$ be an adapted clopen set   for the action $\Phi$, with stabilizer subgroup $G_U \subset G$.   Then for $g, g' \in G$ with $g \cdot U \cap g' \cdot U \ne \emptyset$ we have $g^{-1} \, g' \cdot U = U$, hence $g^{-1} \, g' \in G_U$. Thus,  the  translates $\{ g \cdot U \mid g \in G\}$ form a finite clopen partition of $\fX$, and are in 1-1 correspondence with the quotient space $X_U = G/G_U$ so the stabilizer group $G_U \subset G$ has finite index.  
  Note that the action of $g \in G_U$ on $X_U$ is trivial precisely when $g \in C_U$, where $C_U \subset G_U$ is the largest normal subgroup of $G$ contained in $G_U$. Thus, the action of    the finite group $G/C_U$ on $X_U$ by permutations is a finite approximation of the action of $G$ on $\fX$.
 
The restricted action of $G_U$ on $U$ defines  a homomorphism $\Phi_U \colon G_U \to \Homeo(U)$, with kernel $\ker (\Phi_U) \subset G_U$. Suppose that $g \in \ker (\Phi_U)$ and the action $\Phi$ is quasi-analytic, then $g$ acts trivially on $\fX$, which implies that $g$ acts trivially on $X_U$ hence  $g \in C_U$. Thus, the quasi-analytic hypothesis is equivalent to the assumption that $\ker (\Phi_U)$ is a normal subgroup of $G$.
 Exactly the same reasoning yields the following criteria for an action to be locally quasi-analytic. 
 
 \begin{prop}\label{prop-notLQA}
 Let  $(\fX,G,\Phi)$   be a minimal equicontinuous  action  on a Cantor space $\fX$. Suppose that the action $\Phi$ is locally quasi-analytic. Then there exists $\e > 0$ such that for each pair of non-empty adapted clopen sets   $V \subset U \subset \fX$ with $\diam(U) < \e$, the subgroup  $\ker(\Phi_{V}) \subset   G_V \subset G_U$   is normal in $G_U$.
 \end{prop}
 The action of the Grigorchuk group on the tree boundary (as discussed in Example~\ref{ex-arboreal}) admits     pairs of  clopen subsets $V \subset U$ with arbitrarily small diameters, for which $\ker (\Phi_V)$ is not normal in $G_U$, and thus  its action  is not locally quasi-analytic.

 We next develop the idea behind the proof of Proposition~\ref{prop-notLQA} to develop an effective criteria, given in Theorem~\ref{thm-LQA},  for showing that an action must be   locally quasi-analytic.  
 
 For a choice of basepoint $x \in \fX$ and scale $\e > 0$, there exists an adapted clopen set $U \in \CO(\fX)$ with $x \in U$ and $\diam(U) < \e$.      Iterating this construction, for a given basepoint $x$, one can always construct the following:
 
\begin{defn}\label{def-adaptednbhds}
Let  $(\fX,G,\Phi)$   be a minimal equicontinuous  action  on a Cantor space $\fX$.
A properly descending chain of clopen sets $\cU = \{U_{\ell} \subset \fX  \mid \ell \geq 1\}$ is said to be an \emph{adapted neighborhood basis} at $x \in \fX$ for the action $\Phi$  if
    $x \in U_{\ell +1} \subset U_{\ell}$ for all $ \ell \geq 1$ with     $\cap    U_{\ell} = \{x\}$, and  each $U_{\ell}$ is adapted to the action $\Phi$.
\end{defn}

 For such a collection, set $G_{\ell} = G_{U_{\ell}}$     we obtain a descending chain of finite index subgroups 
 $$\ds \cG_{\cU} = \{G = G_0 \supset G_1 \supset G_2 \supset \cdots \} \ .$$
 The intersection $\ds K(\cG_{\cU}) = \bigcap_{\ell \geq 0} ~ G_{\ell}$ is called the kernel of $\cG_{\cU}$.   

Next, set $X_{\ell} = G/G_{\ell}$ and note that  $G$ acts transitively on the left on   $X_{\ell}$.    
The inclusion $G_{\ell +1} \subset G_{\ell}$ induces a natural $G$-invariant quotient map $p_{\ell +1} \colon X_{\ell +1} \to X_{\ell}$.
 Introduce the inverse limit 
 \begin{equation} \label{eq-invlimspace}
X_{\infty} \equiv \varprojlim \ \{p_{\ell +1} \colon X_{\ell +1} \to X_{\ell} \mid \ell > 0\} 
\end{equation}
which is a Cantor space with the Tychonoff topology, and the action on the factors $X_{\ell}$ induces   a minimal  equicontinuous action $\Phi_x \colon G \times X_{\infty} \to X_{\infty}$. 
 Note that for $g \in K(\cG_{\cU})$, the left action of $g$ on $X_{\ell}$ fixes the coset $e_{\ell} \in X_{\ell}$ and hence fixes the limiting point $e_{\infty} \in X_{\infty}$.

 For each $\ell \geq 0$, we have the ``partition coding map'' $\sigma_{\ell} \colon \fX \to X_{\ell}$ which is $G$-equivariant. The maps $\{\sigma_{\ell}\}$ are compatible with the quotient maps in \eqref{eq-invlimspace}, and so define a  limit map $\sigma_{\infty} \colon \fX \to X_{\infty}$. The fact that the diameters of the clopen sets $\{U_{\ell}\}$ tend to zero, implies that $\sigma_{\infty}$ is a homeomorphism. Let $\tau_x \colon X_{\infty} \to \fX$ denote the inverse map, which commutes with the $G$ actions on both spaces, and satisfies $\tau_x(e_{\infty}) = x$.  The minimal equicontinuous action $(X_{\infty}, G, \Phi_x)$   is called the \emph{odometer representation} centered at $x$ for the action $(\fX,G,\Phi)$.
 
Suppose that $\cU$ and $\cU'$ are two choices of adapted neighborhood bases at $x$, then the corresponding   group chains  $\cG_{\cU}$ and $\cG_{\cU'}$  are equivalent as descending chains, as shown in \cite{FO2002,DHL2016a}. For two choices of basepoints $x,x' \in \fX$, with adapted neighborhood basis $\cU$ at $x$ and $\cU'$ at $x'$,   
the corresponding group chains $\cG_{\cU}$ for $\cU$ and $\cG_{\cU'}$ for $\cU'$ need not be equivalent, but  are \emph{conjugate equivalent} as descending chains.  The conjugacy relation on group chains, as introduced by Fokkink and Oversteegen in \cite{FO2002}, is in essence  the relation on group chains which results in conjugating them by elements of the closure of the group action. It   is key for showing in \cite{DHL2016a} that the invariants  of the action $(\fX,G,\Phi)$ defined in \cite{DHL2016c,HL2017a} are independent of the choices of basepoints and group chain models.

We recall a basic property of the odometer models:
\begin{lemma}
A minimal equicontinuous Cantor action $(\fX,G,\Phi)$ is   free   if and only if for all $x \in \fX$ and  any adapted neighborhood basis $\cU$ at $x$, the kernel $K(\cG_{\cU})$  is the trivial group. 
\end{lemma}
Let $\cG_U$ be the group chain associated to adapted neighborhood system at $x$, and suppose that     $y = g \cdot x$, then set 
$$\cG_U^y =   \{G = G_0 \supset g \, G_1 \, g^{-1} \supset g \, G_2 \, g^{-1}  \supset \cdots \} \ , $$
which is the group chain associated to the adapted neighborhood basis $\cU^g$ at $y$. Clearly, $K(\cG_{\cU}^y) = g \, K(\cG_{\cU}^x) \, g^{-1}$, so the property that $K(\cG_{\cU}^y)$ is trivial   is independent of the choice of $y \in \cO(x)$. On the other hand, if $y \not\in \cO(x)$ then whether $K(\cG_{\cU}^y)$ is trivial or not may depend on the choice of $y$.

We next prove a criterion for when a   minimal equicontinuous Cantor action $(\fX,G,\Phi)$   is locally quasi-analytic.
 Recall the following property of groups.
\begin{defn}\cite{Baer1956}\label{def-noetherian} 
A group $\G$ is said to be \emph{Noetherian} if every increasing chain of closed subgroups $\{H_i \mid i \geq 1 \}$ of $\G$ has a maximal element $H_N$.
\end{defn}
An equivalent definition of the Noetherian property is that every subgroup of   $\G$ is finitely generated. Hence, $\G$ has  at most countably many subgroups. 

  Recall that a group $\G$ is \emph{polycyclic} if there exists an integer $k > 0$ and a chain of subgroups
\begin{equation}\label{eq-polycyclic}
\{e\} = \G_{k+1} \subset \G_{k} \subset \cdots \subset \G_0  = \G
\end{equation}
such that each $\G_{\ell +1}$ is normal in $\G_{\ell}$ and the quotient $\G_{\ell}/\G_{\ell +1}$ is a cyclic group for $0 \leq \ell < k$.
For example, a finitely-generated nilpotent group is polycyclic. 
A group $\G$ is \emph{virtually polycyclic} if there exists a subgroup $\G_0 \subset \G$ of finite index   such that  $\G_0$ is polycyclic.  The following result is folklore.
\begin{prop}\label{prop-polycyclic}
Let $\G$ be a virtually polycyclic group, then $\G$ is Noetherian.
\end{prop}

The following result is a consequence of the Noetherian property for Cantor actions.
 
  \begin{thm}\label{thm-LQA}
Let    $G$ be a Noetherian group. Then   a  minimal equicontinuous Cantor action $(\fX,G,\Phi)$   is locally quasi-analytic.
\end{thm}
 \proof
We assume that the action $(\fX,G,\Phi)$   is not locally quasi-analytic, and  construct an increasing chain of subgroups in $G$ with no maximal element,   contradicting that $G$ is   Noetherian.

 Fix $x \in \fX$ and let $\cU = \{U_{\ell} \subset \fX  \mid \ell \geq 1\}$  be an  adapted neighborhood basis  at $x$ for the action $\Phi$. For $\ell \geq 1$,   the collection of translates $\{g \cdot U_{\ell} \mid g \in G\}$ is a disjoint clopen covering of $\fX$. Let $\lambda_{\ell} > 0$ be a Lebesgue number for this covering. That is, if $U \subset \fX$ is an open set with $\diam (U) < \lambda_{\ell}$ then there exists $g \in G$ such that $U \subset g \cdot U_{\ell}$.

The assumption that Definition~\ref{def-LQA} does not hold  implies that for each $\e > 0$,  there exists  non-empty clopen subsets $V \subset U \subset \fX$ with $\diam(U) < \e$ such that 
\eqref{eq-lqa} fails for this pair. That is, there exists $h_1, h_2 \in G$ such that $\Phi(h_1)|V = \Phi(h_2)|V$ but  $\Phi(h_1)|U \ne \Phi(h_2)|U$. Set $h = h_2^{-1} \, h_1$ then $\Phi(h) |V$ is   the identity map, but $\Phi(h) | U$ is not. Note that for this $h$ we have $h \cdot U \cap U \ne \emptyset$.

 We next construct an increasing chain of subgroups $\{H_i \mid i \geq 1\}$ by induction. 
To begin, set   $\ell_1 =1$.
Let $W_1 \subset \fX$ be a non-empty open set with $\diam (W_1) < \lambda_1$ such that there exists a non-empty open subset $V_1 \subset W_1$ and $h_1 \in G$ so that  $\Phi(h_1) | V_1$ is the identity, and $\Phi(h_1) | W_1$ is not, with $h_1 \cdot W_1 \cap W_1 \ne \emptyset$.
 Then there exists $g_1 \in G$ such that $W_1 \subset g_1 \cdot  U_1$, so $g_1^{-1} \cdot W_1 \subset U_1$. Set $\xi_1 = g_1^{-1} \, h_1 \, g_1$ and note that $\xi_1 \cdot U_1 \cap U_1 \ne \emptyset$, hence $\xi_1 \in G_{U_1}$.   
   
 Note that  $g_1^{-1} \cdot V_1 \subset U_1$. As the action of $G_{U_1}$ on $U_1$ is minimal, there exists $\gamma_1 \in G_{U_1}$ such that $\gamma_1 \cdot x \in   g_1^{-1} \cdot V_1$. Set $V_1' = \gamma_1^{-1} \,  g_1^{-1} \cdot V_1$ so that $x \in V_1' \subset U_1$. 
 Then for 
 $$k_1 =     \gamma_1^{-1} \, \xi_1 \, \gamma_1 =  \gamma_1^{-1} \, g_1^{-1} \, h_1 \, g_1 \, \gamma_1 \ ,$$
 we have that $\Phi(k_1) | U_1$ is not the identity, but $\Phi(k_1) | V_1'$ is the identity.
 Choose $\ell_2 > \ell_1 = 1$ such that $U_{\ell_2} \subset V_1'$ and thus $\Phi(k_1) | U_{\ell_2}$ is also the identity.
 
 Assume that for $i > 1$, and increasing sequence of integers $1 = \ell_1 < \ell_2 < \cdots < \ell_i$   have been chosen as above. 
Let $W_i \subset \fX$ be a non-empty open set with $\diam (W_i) < \lambda_{\ell_i}$ such that there exists a non-empty open subset $V_i \subset W_i$ and $h_i \in G$ so that  $\Phi(h_i) | V_i$ is the identity, and $\Phi(h_i) | W_i$ is not the identity, with $h_i \cdot W_i \cap W_i \ne \emptyset$.
Then there exists $g_i \in G$ such that $W_i \subset g_i \cdot  U_{\ell_i}$ and so $g_i^{-1} \cdot W_i \subset U_{\ell_i}$. 
Set $\xi_i = g_i^{-1} \, h_i \, g_i$ and note that $\xi_i \cdot U_{\ell_i} \cap U_{\ell_i} \ne \emptyset$, hence $\xi_i \in G_{U_{\ell_i}}$.  
 
  Note that  $g_i^{-1} \cdot V_i \subset U_{\ell_i}$. 
   As the action of $G_{U_{\ell_i}}$ on $U_{\ell_i}$ is minimal, there exists $\gamma_i \in G_{U_{\ell_i}}$ such that $\gamma_i \cdot x \in   g_i^{-1} \cdot V_i$. 
     Set $V_i' = \gamma_i^{-1} \,  g_i^{-1} \cdot V_i$ so that $x \in V_i' \subset U_{\ell_i}$. 
  Then for 
 $$k_i =     \gamma_i^{-1} \, \xi_i \, \gamma_i =  \gamma_i^{-1} \, g_i^{-1} \, h_i \, g_i \, \gamma_i \ ,$$
 we have that $\Phi(k_i) | U_{\ell_i}$ is not the identity, but $\Phi(k_i) | V_i'$ is the identity.
   Choose $\ell_{i+1} > \ell_i$ such that $U_{\ell_{i+1}} \subset V_i'$ and thus $\Phi(k_i) | U_{\ell_{i+1}}$ is also the identity.
   
Now assume that indices $\{\ell_i \mid i \geq 1\}$ and elements $\{k_i \mid i \geq 1\}$ have been chose as above. Define:
\begin{equation}
H_i = \{ g \in G_{\ell_i} \mid \Phi(g) | U_{\ell_i} = id  \} \ .
\end{equation}
Observe that for each $g \in H_i$ and $j > i$ we have $U_{\ell_{j}} \subset U_{\ell_i}$ hence   $g \cdot U_{\ell_{j}} = U_{\ell_{j}}$ thus $H_i \subset G_{\ell_{j}}$.
Moreover, $\Phi(g) | U_{\ell_{j}}$ is the identity, hence $H_i \subset H_{j}$. 

Now let $i \geq 1$, then by choice, $\Phi(h_i) | W_i$ is not the identity, hence $\Phi(k_i) | U_{\ell_i}$ is not the identity, so $k_i \not\in H_{i}$. On the other hand, for $j > i$ the restriction   $\Phi(k_i) | U_{\ell_j}$ is   the identity, hence $k_i \in H_j$.
 
It follows that the collection $\{H_i \mid i \geq 1\}$ forms a strictly increasing chain of subgroups in $G$, which therefore is not Noetherian.
\endproof

Theorem~\ref{thm-LQA} implies the following extension of       Corollary~\ref{cor-abelian}.

 \begin{cor}\label{cor-nilLQA}
 Let    $G$ be a finitely generated   nilpotent group. Then   a  minimal equicontinuous Cantor action $(\fX,G,\Phi)$   is locally quasi-analytic.
 \end{cor}

  The proof of Theorem~\ref{thm-LQA} yields the following   result.
\begin{cor}\label{cor-wild}
Let  $(\fX,G,\Phi)$  be a  minimal equicontinuous Cantor action which  is not locally quasi-analytic. Then for any  $x \in \fX$,  the isotropy subgroup $G_x \subset G$ contains an infinite strictly increasing chain of subgroups.
\end{cor}
\proof
For given $x$, choose an adapted neighborhood basis at $x$. Then   for $H_i$ as constructed in the proof above, the action of each $h \in H_i$ fixes the set $U_{\ell_{i+1}}$ so in particular $h \in G_x$. Thus $H_i \subset G_x$.
\endproof

 It is remarkable that the first construction of a non-homogenous weak solenoid by Schori in \cite{Schori1966}, see also \cite{CFL2014}, has this ascending chain property for the isotropy subgroups of its monodromy action. 
 The Grigorchuk groups \cite{Grigorchuk1984}, and more generally branch groups \cite{BGS2003,Grigorchuk2000}, provide a large class of examples of groups acting on trees which give rise to Cantor actions that are not locally quasi-analytic. This is discussed further in  Example~\ref{ex-arboreal}.

The authors, in joint work with Jessica Dyer, gave in   \cite[Example~8.5]{DHL2016a}  examples of minimal equicontinuous Cantor actions $(\fX, G, \Phi)$ where $G$  is a   torsion-free 2-step nilpotent group,   its discriminant  group $\cD_x$ is a  Cantor space, and   the action is   locally quasi-analytic. It was asked in \cite{DHL2016c} whether it is possible to construct  examples of minimal equicontinuous Cantor actions $(\fX, G, \Phi)$ where $G$  is a   finitely-generated nilpotent group, and the action is not locally quasi-analytic. Corollary~\ref{cor-nilLQA}  shows this is impossible.

\section{Continuous orbit equivalence and rigidity} \label{sec-coe}

The concept of  \emph{continuous orbit equivalence} between Cantor actions was introduced by Mike  Boyle in his thesis \cite{Boyle1983}, and has played a fundamental role in  the classification of Cantor actions in many subsequent works \cite{BezuglyiMedynets2008,GPS1999,GMPS2010}. The related notion of the topological full group of an action in Definition~\ref{def-TFG} has   provided a rich source of examples of finitely generated groups with exceptional properties, as discussed for example in \cite{deCornulier2014}. In this section, we   show how these notions can be used to show return equivalence for minimal equicontinuous Cantor actions.

 \begin{defn}\label{def-coe}
For $i=1,2$, let $(\fX_i,G_i,\Phi_i)$   be an action  on a Cantor space $\fX_i$.  The    actions are said to be \emph{continuously orbit equivalent} if there exists a homeomorphism $h \colon \fX_1 \to \fX_2$ and continuous functions
\begin{enumerate}
\item $\alpha \colon G_1 \times \fX_1 \to G_2$,  $h(\Phi_1(g_1)(x_1)) = \Phi_2(\alpha(g_1 , x_1),  h(x_1))$ for all   $g_1 \in G_1$ and $x_1 \in \fX_1$; 
\item $\beta \colon  G_2 \times \fX_2 \to G_1$,  $h^{-1}(\Phi_2(g_2, x_2)) = \Phi_1(\beta(g_2, x_2), h^{-1}(x_2))$ for all $g_2 \in G_2$ and $x_2 \in \fX_2$. 
\end{enumerate}
The    actions are said to be \emph{virtually continuously orbit equivalent} if there exists adapted clopen subsets $U_i \subset \fX_i$ such that the restricted actions $(U_i,G_{U_i},\Phi_i)$ are continuously orbit equivalent.
\end{defn}
 The homeomorphism  $h$ is called a  \emph{continuous orbit equivalence} between the two actions.
Note that  the functions $\alpha$ and $\beta$ are not assumed to satisfy the cocycle property.

   Given an  action  $(\fX_2,G_2,\Phi_2)$ and a   homeomorphism  $h \colon \fX_1 \to \fX_2$ ,   define the conjugate action $(\fX_1,G_2,\Phi_2^h)$ by setting
 $$\Phi_2^h(g_2, x_1) = h^{-1}(\Phi_2(g_2, h(x_1))) ~  {\rm for} ~ g_2 \in G_2, ~ x_1 \in \fX_1 \ . $$
We recall an observation from \cite[Section~2]{GW1995} about the topological full group:
\begin{prop}\label{prop-coetfg}
For $i=1,2$, let $(\fX_i,G_i,\Phi_i)$   be an action  on a Cantor space $\fX_i$. Then a homeomorphism $h \colon \fX_1 \to \fX_2$ is a continuous orbit equivalence between the actions, if and only if the following two conditions hold:
\begin{enumerate}
\item $\Phi_1(g_1, \cdot) \in [[\fX_1,G_2,\Phi_2^h]]$ for all $g_1 \in G_1$;
\item $\Phi_2^h(g_2, \cdot) \in [[\fX_1,G_1,\Phi_1]]$ for all $g_2 \in G_2$.
\end{enumerate}
That is, the actions are continuously orbit equivalent if and only if the images $\Phi_1(G_1) \subset \Homeo(\fX_1)$ and $\Phi_2^h(G_2) = h^{-1} \circ \Phi_2(G_2) \circ h \subset \Homeo(\fX_1)$ generate the same topological full group.
\end{prop}

In the works \cite{CortezMedynets2016,Li2015}, the following    is called  the ``rigidity problem''   for Cantor actions. 

 \begin{prob}
For $i=1,2$, let $(\fX,G_i,\Phi_i)$   be an action  on a fixed Cantor space $\fX$, and suppose that $[[\fX,G_1,\Phi_1]] = [[\fX,G_2,\Phi_2]]$. Give conditions on the actions $\Phi_1$ and $\Phi_2$ which are sufficient to imply that they are return equivalent.
 \end{prob}
 
 Given a Cantor action $(\fX,G_1,\Phi_1)$, one can   construct a new Cantor action $(\fX,G_2,\Phi_2)$ with the same topological full group   by simply adjoining elements of $[[\fX_1,G_1,\Phi_1]]$ to the action $\Phi_1(G_1)$. We   discuss this construction  further  in Example~\ref{ex-modifications}. The construction of such examples shows that   for a solution of the rigidity problem, there must be   hypotheses   imposed   which rule them out. 
  
Let $(\fX_i,G_i,\Phi_i)$   be a minimal equicontinuous  action  on a Cantor space $\fX_i$ for    $i=1,2$.  Assume there exists a  continuous orbit equivalence $h \colon \fX_1 \to \fX_2$ between   two actions, then by the above remarks, we can assume without loss of generality that $\fX_1 = \fX_2 = \fX$,  and  there exists functions $\alpha, \beta$ as in Definition~\ref{def-coe}. If the both actions are free, then the  functions $\alpha$ and $\beta$ are uniquely determined by the actions and satisfy a cocycle property. This is a key point in the works     \cite{CortezMedynets2016,GPS2017}. It was observed in \cite{BoyleTomiyama1998,Li2015}, and also implicitly in \cite{Renault2008}, that if  the actions are topologically free on $\fX$, then the orbit functions $\alpha$ and $\beta$  still   satisfy the cocycle identities. Li showed that if, in addition,  the actions satisfy the  \emph{continuous cocycle rigidity} property of \cite[Section~4]{Li2015}, then the actions are conjugate. 

The work by Cortez and Medynets \cite{CortezMedynets2016} showed that for \emph{free} minimal equicontinuous Cantor actions, the orbit cocycles $\alpha$ and $\beta$ induce a return equivalence between the two actions. 
The   observation behind the proof of Theorem~\ref{thm-main1} is that if  both actions are assumed to be locally quasi-analytic, then the  method of proof for     \cite[Theorem~3.3]{CortezMedynets2016} can be  used to show   the   actions are return   equivalent.

   \proof[Proof of Theorem~\ref{thm-main1}]
Medynets states in \cite[Remark~3]{Medynets2011} (and also in Theorem~4.2 in \cite{CortezMedynets2016}) that if the topological full groups of two minimal Cantor actions  $(\fX,G,\Phi)$  and  $(\fX,H,\Psi)$ are isomorphic as abstract groups, then there is a homeomorphism $h \colon \fX \to \fX$ which induces the isomorphism. That is, every such abstract group isomorphism is implemented by a spatial isomorphism, which by Proposition~\ref{prop-coetfg} yields a topological orbit equivalence between the actions.

Thus,  we may assume   that     $(\fX,G,\Phi)$  and  $(\fX,H,\Psi)$   are   minimal equicontinuous  actions on a Cantor space $\fX$ which have the same orbits. Let   $\e > 0$ be chosen so that  both actions are locally quasi-analytic for clopen sets $U \subset \fX$ satisfying $\diam(U) < \e$.
Choose a basepoint $x \in \fX$. 

Let   $V  \subset \fX$ be an adapted clopen set for the action $\Psi$ with $x \in V$ and $\diam (V) < \e$, with stabilizer group $H_V \subset H$. Then    by Proposition~\ref{prop-BCT},   the action of $\cH^{\Psi}_{V} = \Psi_V(H_V) \subset \Homeo (V)$  on $V$ is   topologically free. 

Let $\cU = \{U_{\ell} \subset \fX  \mid \ell \geq 1\}$ 
be an  adapted neighborhood basis at  $x$  for the action  $(\fX,G,\Phi)$, with  $U_1\subset V$.  
 Let $\cG_{\cU} = \{G_{\ell} \mid \ell \geq 0\}$ be the associated group chain with $G_0 = G$ and $G_1 = G_{U_1}$.
   Then  the action of $\cH^{\Phi}_{U_1} = \Phi_{U_1}(G_1) \subset \Homeo (U_1)$  on $U_1$ is also  topologically free.

Let $\alpha \colon  G  \times \fX  \to H$ be the continuous   function in Definition~\ref{def-coe}.
The subgroup $G_1 \subset G$ has finite index, and $G$ is finitely-generated, so there exists   a finite generating set $\{g_1, \ldots , g_k\} \subset G_1$. 
By the continuity of $\alpha$, there exists  $\delta_1 > 0$ such that
\begin{equation}\label{eq-locconstant}
  \alpha(g_i, x) = \alpha(g_i, y) ~ {\rm    for } ~    1 \leq i \leq k, ~ x,y \in \fX ~ {\rm with} ~ \dX(x,y) < \delta_1 \ .
\end{equation}
 
 Let $\ell_2 \geq 2$ be such that $\diam (g \cdot U_{\ell_2}) < \delta_1$ for all $g \in G_1$. The collection $\{g \cdot U_{\ell_2} \mid g \in G_1\}$ is a finite clopen partition of $U_1$ so there exists $0 < \delta_2 < \delta_1$ such that for any $y \in U_1$   there exists $g \in G_1$ such that $B_{\dX}(y,\delta_2) \subset g \cdot U_{\ell_2}$. Then $\delta_2$ is a Lebesgue number for the open covering of $U_1$ by translates of $U_{\ell_2}$.

By the uniform continuity of the action $\Phi$, there exists  $0 < \delta_3 \leq \delta_2$ such that, for all $g \in G$, 
\begin{equation}
\dX(\Phi(g)(x) ,\Phi(g)(y)) < \delta_2 ~ {\rm for ~ all} ~  x,y \in \fX ~ {\rm with} ~  \dX(x,y) < \delta_3 \ .
\end{equation}
  
  Let $\ell_3 > \ell_2$ be such that $\diam (g \cdot U_{\ell_3}) < \delta_3$ for all $g \in G_1$.  Note that $U_{\ell_3} \subset U_{\ell_2} \subset U_1$.
    
 Now consider the restriction $\alpha \colon  G_1 \times U_{1} \to H$. For each $g \in G_1$ and $x \in U_{1}$ we have $\Phi(g)(x) \in U_1$. Let $h \in H$ be such that $\Psi(h)(x) = \Phi(g)(x)$, then $U_1 \subset V$ implies that $\Psi(h)(V) \cap V \ne \emptyset$ hence $h \in H_V$.  
 Recall that $\cH^{\Psi}_{V} = \Psi_V(H_V) \subset \Homeo(V)$,  thus 
  the restriction of $\alpha$ to $G_1 \times U_{1}$ induces a map  $\whalpha_1 = \Psi_V \circ \alpha \colon  G_1 \times U_{1} \to \cH^{\Psi}_{V}$.

 \begin{lemma}\label{lem-cocycle}
 The map $\whalpha_1 \colon  G_1 \times U_1 \to \cH^{\Psi}_{V}$ is a cocycle over the action of $\Phi_{U_1}$ on $U_1$.
 \end{lemma}
 \proof
  Let $g,g' \in G_1$ and let $x \in  U_1$. Then $\Phi_{U_1}(g' \, g)(x) = \Phi_{U_1}(g') (x')$ for $x' = \Phi_{U_1}(g)(x)$.
  
  The  action of $\cH^{\Psi}_{V}$ is topologically free, so there exists 
 a dense $\cH^{\Psi}_{V}$-invariant subset $Z \subset V$ such that the action of $\cH^{\Psi}_{V}$ on $Z$ is free. 
  That is, for $\psi, \psi' \in  \cH^{\Psi}_{V}$ and $z \in Z$, if $\psi(z) = \psi'(z)$ then    $\psi = \psi' \in \Homeo(V)$. 
  
  Let $x \in Z$ and $h= \alpha_1(g,x)$, then $\psi = \Psi_V(h)$ is the unique map in $\cH^{\Psi}_{V}$ such that $\psi(x) = x'$. Note that $x' \in Z$ so there is a unique $\psi' \in \cH^{\Psi}_{V}$ such that $x'' = \psi'(x') = \Phi_{U_1}(g')(x')$. Then $x'' = \psi' \circ \psi(x)$.

By  the defining identity (1) in   Definition~\ref{def-coe} with $h$ the identity map, for $x \in Z \cap U_1$ we have that 
 $$\whalpha_1(g , x) = \psi ~ , ~  \whalpha_1(g' , x') = \psi'  ~ , ~ \whalpha_1(g' \, g , x) = \psi' \,  \psi  \ .$$
 The continuity of the actions $\Phi$ and $\Psi$ imply that    this holds for all $x \in  U_1$ as $Z \cap U_1$ is dense in $U_1$.   
 Thus,  we have the cocycle identity 
  \begin{equation}
\whalpha_1(g' \, g , x) = \whalpha_1(g'  , \Phi_{U_1}(x)) \, \whalpha_1(g , x) 
\end{equation}
for all $g', g \in G_{U_1}$ and $x \in U_1$, as was to be shown.
\endproof

  Next, we show that the cocycle $\widehat{\alpha}_1$ is defined by a group homomorphism     when restricted to the clopen subset $U_{\ell_3} \subset U_1$. 
  The method of proof is   the same as for the proof of \cite[Theorem~3.3]{CortezMedynets2016}.

Recall that $\{g_1, \ldots , g_k\} \subset G_1$ is a   set of generators, then 
  for each $1 \leq i \leq k$,   by the choice of $U_{\ell_2}$ and \eqref{eq-locconstant},  for any $g' \in G_1$ the value of $\whalpha_1(g_i,x)$ is constant for $x \in g' \cdot U_{\ell_2}$. 

By the choice of $\delta_1$,   for $\dX(x,y) < \delta_1$ and $1 \leq i \leq k$,   we have   $\widehat{\alpha}_1(g_i ,x) = \widehat{\alpha}_1(g_i ,y)$.
Then by the choice of $\delta_3$ and $U_{\ell_3}$,  for any $1 \leq i \leq k$,  $g' \in  G_{U_1}$ and $x,y \in  U_{\ell_3}$ then
$$ \widehat{\alpha}_1(g_i , \Phi_{U_1}(g')(x)) = \widehat{\alpha}_1(g_i , \Phi_{U_1}(g')(y)) \ .$$
We apply this to the case where $g = g_{i_1} \cdots g_{i_{\nu}} \in G_{U_{\ell_3}}$. 
Set $g' = g_{i_2} \cdots g_{i_{\nu}}$ then for $x,y \in  Z \cap U_{\ell_3}$, 

\begin{eqnarray}
 \lefteqn{ \widehat{\alpha}_1(g , x) = \widehat{\alpha}_1(g_{i_1} \,  g', x)   =    \widehat{\alpha}_1(g_{i_1} , \Phi_{U_1}(g')(x)) \circ  \widehat{\alpha}_1(g' , x) } \label{eq-computation}\\
& = &       \widehat{\alpha}_1(g_{i_1} , \Phi_{U_1}(g')(y)) \circ  \widehat{\alpha}_1(g' , y) 
 =   \widehat{\alpha}_1(g_{i_1} \,  g', y)
 = \widehat{\alpha}_1(g , y) . \nonumber
\end{eqnarray}

As the values of $\widehat{\alpha}_1(g , x), \widehat{\alpha}_1(g , y) \in \cH^{\Psi}_{V}$ are defined using the identity (1) in   Definition~\ref{def-coe}, the identity \eqref{eq-computation} holds   on the closure of $Z \cap U_{\ell_3}$ which is all of $U_{\ell_3}$.
Thus, the     calculation \eqref{eq-computation} shows that the restricted cocycle 
$\ds \whalpha_3 \colon G_{U_{\ell_3}} \times U_{\ell_3} \to \cH^{\Psi}_{V}$ 
is independent of the point in the second factor $U_{\ell_3}$,  hence induces a   group homomorphism denoted by 
$\ds  \wtalpha_3 \colon G_{U_{\ell_3}} \to \cH^{\Psi}_{V}$. 

 Suppose  that   $g, g' \in G_{U_{\ell_3}}$ satisfy $\Phi_{U_{\ell_3}}(g) = \Phi_{U_{\ell_3}}(g') \in \Homeo(U_{\ell_3})$, then by  the defining identity (1) in   Definition~\ref{def-coe} with $h$ the identity map, for $x \in Z \cap U_1$ we have  $\whalpha_1(g,x) = \whalpha_1(g',x)\in  \cH^{\Psi}_{V}$.
 It follows that $\whalpha_3$ induces a homomorphism 
$\whalpha_3 \colon  \cH^{\Phi}_{U_{\ell_3}}  \to \cH^{\Psi}_{V}$. 
The homomorphism $\whalpha_3$ is in fact an isomorphism.  The proof is  omitted, as it follows by   the same arguments as in the works \cite{Li2015,CortezMedynets2016}.

Set ${\bf G} = \cH^{\Phi}_{U_{\ell_3}}$ and ${\bf H} = \whalpha_3({\bf G}) \subset \cH^{\Psi}_{V}$, and denote by ${\bf A} \equiv \whalpha_3 \colon {\bf G} \to{\bf H}$. 

Set $W =  U_{\ell_3}$.  Note that the group ${\bf G}$ is the stabilizer of  $W$, and the orbits of ${\bf H}$ on points in $W$ equal the orbits of ${\bf G}$, hence ${\bf H}$ also stabilizes $W$.
Thus,  ${\bf A}$ conjugates the   action of $\bf G$  and $\bf H$ on  $W$.
This completes the proof of Theorem~\ref{thm-main1}.
 \endproof

 The conclusion of Theorem~\ref{thm-main1} cannot be improved to obtain a conjugacy of the actions $\Phi$ and $\Psi$, as illustrated in  Examples~\ref{ex-direct} and \ref{ex-semidirect}. However, one can show that $\bf G$ and $\bf H$ determine subgroups of the same index in $G$ and $H$, respectively, using the same methods as in the proof of \cite[Theorem~3.3]{CortezMedynets2016}, so  the actions $\Phi$ and $\Psi$ are structurally conjugate in the sense of \cite{CortezMedynets2016}.
  
  Recall that in Definition~\ref{def-coe}, the notion of virtually continuously orbit equivalent actions was introduced, defined in terms of the induced actions on adapted clopen sets for the actions. There is a corresponding extension to the full groups for a Cantor action:
  \begin{defn}\label{def-virtualfullgroup}
  Let   $(\fX,G,\Phi)$  and  $(\fY,H,\Psi)$ be minimal Cantor actions. Their full groups are said to be \emph{virtually isomorphic} if there   exists adapted clopen subsets $U \subset \fX$ and $V \subset \fY$ such that their topological full groups $[[U,G_U,\Phi_U]]$  and $[[V,H_V,\Psi_V]]$ of their restricted actions are isomorphic, and we write  $[[[\fX,G,\Phi]]] \sim [[[\fY,H,\Psi]]]$.
   \end{defn}
   Note that if $\phi \in [[U,G_U,\Phi_U]]$ then as $\fX - U$ is a clopen set, we can extend the map $\phi$ as the identity on $\fX - U$ to obtain $\wtphi \in [[\fX,G,\Phi]]$. This yields an inclusion  $[[U,G_U,\Phi_U]] \subset [[\fX,G,\Phi]]$ of groups. The image subgroup will not have finite index unless $U = \fX$, so the   terminology ``virtual'' is an abuse of the usual sense  of this terminology, but is chosen as the subgroup $G_U \subset G$ does have finite index.

As for the case of the   topological full groups of   actions,   \cite[Remark~3]{Medynets2011} implies that if $[[U,G_U,\Phi_U]]$ and $ [[V,H_V,\Psi_V]]$ are isomorphic as groups, then there  is a homeomorphism $h \colon U \to V$ which induces the isomorphism of the full groups. Thus, if the actions $(\fX,G,\Phi)$  and  $(\fY,H,\Psi)$ have virtually isomorphic topological full groups, then  the restricted actions $(U,G_U,\Phi_U)$  and $(V,H_V,\Psi_V)$ are topologically orbit equivalent. Then by 
 Theorem~\ref{thm-main1} applied to these restricted actions, they are return equivalent. Moreover, if each action is locally quasi-analytic, then the restricted actions are as well. Furthermore, if $W \subset U$ is an adapted set for the restricted action $\Phi_U$ then it is also adapted for the action $\Phi$, and similarly for the action $\Psi_V$.
We thus obtain:
 \begin{thm}\label{thm-VFHrigidity}
 Let   $(\fX,G,\Phi)$  and  $(\fY,H,\Psi)$ be minimal Cantor actions which are locally quasi-analytic, and such that  $[[[\fX,G,\Phi]]] \sim [[[\fY,H,\Psi]]]$. Then the actions are return equivalent.
 \end{thm}
 
\begin{remark}\label{rmk-atfg} {\rm 
  The virtual isomorphism class $[[[\fX,G,\Phi]]]$ of the full group of a minimal Cantor action $(\fX,G,\Phi)$ is a ``local invariant'' of the action, as it is determined by the action restricted to arbitrarily small adapted sets for the action. Another such local invariant    is the asymptotic discriminant   introduced   in \cite{HL2017a},   which distinguishes  large classes of minimal equicontinuous actions.  It  seems an interesting problem  to determine the ``asymptotic invariants'' of return equivalence directly from group invariants of the virtual isomorphism class $[[[\fX,G,\Phi]]]$.   
 }
 \end{remark}

  \section{Classification of nil-solenoids} \label{sec-nil}

 In this section, we   recall the construction and a few basic  properties of  weak solenoids, which are a special class continua    introduced by McCord \cite{McCord1965}. Schori   gave an example in \cite{Schori1966}  of a weak solenoid whose Cantor fiber was not a Cantor group,   which motivated  the study   of the self-homeomorphism group of weak solenoids. This   was further investigated by Fokkink and Oversteegen in \cite{FO2002}, who introduced the methods that were further developed  in the works by the authors   \cite{DHL2016a,DHL2016b,DHL2016c}. 
 
The classification of weak solenoids was studied in the works   \cite{AartsFokkink1991,CHL2017,HL2017a}, where it is shown in particular  that a homeomorphism between weak solenoids induces a return equivalence between their global monodromy Cantor actions. This result is  recalled as Theorem~\ref{thm-re} below.  
 We    introduce the  class of nil-solenoids, then give the proof of Theorem~\ref{thm-main3},   which   represents a broad generalization of the classification result for $1$-dimensional solenoids by  Aarts and Fokkink     in \cite{AartsFokkink1991}.
 
 We first recall some basic concepts of weak solenoids.
Let $M_0$ be a closed connected manifold, and $x_0 \in M_0$ a choice of basepoint. Let $G = \pi_1(M_0, x_0)$ denote its fundamental group. Let $\cG$  be a  properly descending chain of finite index subgroups, 
\begin{equation}\label{eq-groupchain}
\cG = \{G = G_0 \supset G_1 \supset G_2 \supset \cdots \} \ .
\end{equation}
    For $\ell \geq 0$, each subgroup $G_{\ell}$ determines a finite covering $\pi_{\ell} \colon M_{\ell} \to M_0$, where $M_{\ell}$ is a closed manifold.  The inclusions $G_{\ell +1} \subset G_{\ell}$ induce non-trivial proper covering maps $p_{\ell+1} \colon M_{\ell+1} \to M_{\ell}$. The collection of these maps,  $\cP = \{ p_{\ell+1} \colon M_{\ell+1} \to M_{\ell} \mid \ell \geq 0\}$, is called a \emph{presentation}. 

     Associated to  $\cP$ is the \emph{weak solenoid} $\cS_{\cP}$ which is the inverse limit space   
\begin{equation}\label{eq-presentationinvlim}
\cS_{\cP} \equiv \lim_{\longleftarrow} ~ \{ p_{\ell +1} \colon M_{\ell +1} \to M_{\ell}\} ~ \subset \prod_{\ell \geq 0} ~ M_{\ell} ~ .
\end{equation}
 By definition, for a sequence $\{x_{\ell} \in M_{\ell} \mid \ell \geq 0\}$, we have 
\begin{equation}\label{eq-presentationinvlim2}
x = (x_0, x_1, \ldots ) \in \cS_{\cP}   ~ \Longleftrightarrow  ~ p_{\ell}(x_{\ell}) =  x_{\ell-1} ~ {\rm for ~ all} ~ \ell \geq 1 ~. 
\end{equation}
The set $\cS_{\cP}$ is given  the relative  topology, induced from the product topology, so that $\cS_{\cP}$ is  compact and connected. The   map $\Pi_{\ell} \colon \cS_{\cP} \to M_{\ell}$ is given by projection onto the $\ell$-th factor in \eqref{eq-presentationinvlim}. 
  
For example, if    $M_{\ell} = \mS^1$ for each $\ell \geq 0$, and the map $p_{\ell}$ is a proper covering map of degree $m_{\ell} > 1$ for $\ell \geq 1$, then $\cS_{\cP}$ is an example of a   classic solenoid.
  The generalization of $1$-dimensional solenoids to the class of weak solenoids was introduced   by McCord in \cite{McCord1965}. In particular, McCord showed  that   $\cS_{\cP}$ has a uniform local product structure.

\begin{prop}\cite{McCord1965,ClarkHurder2013} \label{prop-solenoidsMM}
Let   $\cS_{\cP}$ be a weak solenoid, whose    base space $M_0$   is a compact manifold of dimension $n \geq 1$. Then  $\cS_{\cP}$  is a foliated space, with foliation $\F_{\cP}$ whose leaves have dimension $n$. That is, for each  $x \in \cS_{\cP}$ there is an open neighborhood   $x \in U_x \subset \cS_{\cP}$ homeomorphic to the product space $(-1,1)^n \times K_x$ where $K_x$ is a Cantor space, which is a foliation chart  for    $\F_{\cP}$.
\end{prop}
A weak solenoid is a  \emph{matchbox manifold} of dimension $n$ in the terminology of \cite{ClarkHurder2013}, or a \emph{solenoidal manifold} in the terminology of   \cite{Sullivan2014,Verjovsky2014}.   

Let $\fX_0 = \Pi_0^{-1}(x_0)$ denote the fiber over the basepoint $x_0 \in M_0$. The global monodromy of the foliation $\F_{\cP}$ on $\cS_{\cP}$  is the   action $\Phi_0 \colon G \times \fX_0 \to \fX_0$   defined by the holonomy transport along the leaves of the foliation $\F_{\cP}$. This action is minimal and equicontinuous. A choice of basepoint $x = (x_0,x_1, x_2, \ldots) \in \fX_0$ defines basepoints in each covering space $M_{\ell}$.

Let $X_{\infty}$ be the inverse limit sequence associated to the group chain $\cG$ in \eqref{eq-groupchain}, as defined by \eqref{eq-invlimspace},  and let 
$\Phi \colon G \times X_{\infty} \to X_{\infty}$ be the associated minimal equivariant Cantor action. 
Then the homeomorphism  $\tau_x \colon X_{\infty} \to \fX_0$ defined in Section~\ref{sec-stable} conjugates the actions $(\fX_0, G , \Phi_0)$ and $(X_{\infty} , G , \Phi)$. Moreover, the kernel  $\ds K(\cG)$ is isomorphic to the fundamental group $\pi_1(L_x, x)$ of the leaf $L_x \subset \cS_{\cP}$ containing $x$.

   \begin{thm} \cite[Theorem~1.1]{CHL2017} \label{thm-re}
 Let $\cP$ be a presentation of a weak solenoid   $\cS_{\cP}$ over a closed manifold $M_0$ of dimension $n$ with fundamental group $G$, and   $\cP'$ be a presentation of a weak solenoid   $\cS_{\cP'}$ over a closed manifold $M_0'$ of dimension $n'$ with fundamental group $G'$. Suppose that $\cS_{\cP}$ is homeomorphic to $\cS_{\cP'}$, then the global holonomy action $(X_{\infty}, G , \Phi)$ of $\F_{\cP}$ is return equivalent to the global holonomy action $(X_{\infty}', G' , \Phi')$ of $\F_{\cP'}$. 
   \end{thm}
 We introduce the following terminology, which is well-defined by   Theorem~\ref{thm-re}.
 \begin{defn}
Let  $[[[\F_{\cP} ]]]$ denote the virtual isomorphism class   of the full group of the monodromy action $(\fX_0, G , \Phi_0)$ for the weak solenoid  $\cS_{\cP}$ with foliation $\F_{\cP}$.
 \end{defn}
 
    \begin{cor}\label{cor-re}
 Let  $\cS_{\cP}$ and $\cS_{\cP'}$ be homeomorphic weak solenoids, then  $[[[\F_{\cP} ]]] \sim [[[\F_{\cP'}]]]$.
 \end{cor}
 
We can now give an application of Theorem~\ref{thm-main1} which gives a partial converse to Corollary~\ref{cor-re}.
    
   A closed manifold $M_0$ is said to be a \emph{nilmanifold} if its fundamental group $G = \pi_1(M_0, x_0)$ is a nilpotent group, for some choice of basepoint $x_0 \in M_0$, and its universal covering $\wtM_0$ is a contractible space. In particular, this implies that the fundamental group $G$   is  torsion free.     Let $\cG$ be a group chain as in \eqref{eq-groupchain},     such that the kernel  $\ds K(\cG)$ is the trivial group. The inverse limit space 
    $\cS_{\cP}$ is called a \emph{nil-solenoid}.
 
\begin{thm}\label{thm-conjugate}
 Let  $\cS_{\cP}$ and $\cS_{\cP'}$ be nil-solenoids, and suppose that $[[[\F_{\cP} ]]] \sim [[[\F_{\cP'}]]]$.  Then  $\cS_{\cP}$ and $\cS_{\cP'}$  are homeomorphic continua.
\end{thm}
\proof
Let $G_0$ denote the fundamental group for the closed manifold $M_0$ in the presentation $\cP$, and $G_0'$ the fundamental group    for the closed manifold $M_0'$ in the presentation $\cP'$.

The assumption that $[[[\F_{\cP} ]]] \sim [[[\F_{\cP'}]]]$ implies there exists adapted clopen sets $U \subset X_{\infty}$ and $U' \subset X_{\infty}'$ and a homeomorphism $h \colon U \to U'$ that induces an isomorphism between the topological full groups 
$[[U, G_U , \Phi_{U}]]$ and $[[U', G'_{U'} , \Phi'_{U'}]]$ for the minimal equicontinuous Cantor actions $(U, G_U , \Phi_{U})$ and $(U', G'_{U'} , \Phi'_{U'})$.

As $G_0$ and $G_0'$ are torsion-free, nilpotent, and finitely generated groups, the same also holds for the subgroup $G_U$ of finite index, and similarly for $G'_{U'}$.  It then follows from Corollary~\ref{cor-nilLQA} that both restricted actions are locally quasi-analytic. It then follows from Theorem~\ref{thm-main1} that the restricted actions are return equivalent, and thus the actions $(X_{\infty}, G , \Phi)$ and $(X_{\infty}', G' , \Phi')$ are return equivalent.

Note that the assumption in the definition of a nil-solenoid $\cS_{\cP}$ that the foliation $\F_{\cP}$ has a simply connected leaf $L_0 \subset \cS_{\cP}$ implies that the kernel $K(\cG)$ for its defining group chain $\cG$ is trivial. The simply connected leaf $L_0$ is homeomorphic to the universal covering of the base manifold $M_0$ hence is contractible. Analogous comments hold for the nil-solenoid $\cS_{\cP'}$.

Return equivalence implies that the fundamental group $G_0 = \pi_1(M, x_0)$ contains a subgroup of finite-index which is isomorphic to a subgroup of finite index in  $G_0' = \pi_1(M', x_0')$. For a nil-manifold, its dimension is determined by the cohomological dimension of its fundamental group, so $M_0$ and $M_0'$ have finite coverings of the same dimension, hence must have equal dimensions. 

Then by Theorem~1.5 of \cite{CHL2017}, the  continua $\cS_{\cP}$ and $\cS_{\cP'}$ are homeomorphic. 
\endproof

  \appendix
  
 \section{Examples} \label{sec-examples}

In this appendix, we give a collection of examples of   minimal  Cantor actions, which   illustrate the complexities that these actions can exhibit. The first examples   are elementary, then increase in the subtleness of their properties. The presentations of the examples are   brief,   with references to their detailed constructions and study in the cited literature, as the interest is in showing the possibilities.
  \begin{ex}[Abelian actions]\label{ex-free}
  {\rm By Corollary~\ref{cor-abelian}, for $G$   a finitely-generated abelian group, every effective minimal action of $G$ on a Cantor space $\fX$ is   quasi-analytic, hence   topologically free. 
The work of Li in \cite{Li2015} studies the relation between continuous orbit equivalence and conjugation of minimal equicontinuous Cantor actions. Example~3.5 in \cite{Li2015} constructs two   abelian actions  of $\mZ^n$ which are continuously orbit equivalent but not conjugate. By Theorem~1.2 in \cite{CortezMedynets2016}, or Theorem~\ref{thm-main1} above, the actions   must be return equivalent. 
 The results of  Giordano,  Putman and  Skau in \cite{GPS2017}   classifies such actions in terms of their Pontrjagin duals.

}
  \end{ex}
  
    \smallskip
    
  \begin{ex}[Direct products]\label{ex-direct} 
 {\rm 
 
 Consider first a ``toy'' example. Let $G$ and $G'$ be non-isomorphic finite groups of the same order, $N = |G| = |G'| \geq 4$. 
 Let $X = \{1,2, \ldots , N\}$ be a finite set, so obviously not a Cantor space. Choose  an isomorphism $\vp  \colon G \to X$ and define the action $(X, G, \Phi)$ by using $\vp$ to conjugate the left action of $G$ on itself to an action on $X$. Likewise, choose an    isomorphism $\vp' \colon G' \to X$ and define the action $(X, G', \Phi')$ by using $\vp'$ to conjugate the left action of $G'$ on itself to an action on $X$. Note that for both of these actions, there is only one orbit of the action. Then both actions are minimal and equicontinuous, and the identity map $X \to X$ is a continuous orbit equivalence. However, the actions $\Phi$ and $\Phi'$ cannot be conjugate as there is no group isomorphism   $h_* \colon G \to G'$.
 
 This simple example can then be modified to obtain examples which show  that the conclusion of Theorem~\ref{thm-main1} cannot be improved to obtain a conjugacy of the actions. Let $(\fY, H, \Psi)$ be any minimal equicontinuous Cantor action of a free abelian group $H$. Set $\fX = X \times \fY$ which is   a Cantor space. Then the product actions $(\fX, G \times H, \Phi \times \Psi)$ and $(\fX, G' \times H, \Phi' \times \Psi)$ are locally quasi-analytic and continuously orbit equivalent, so satisfy the hypotheses of Theorem~\ref{thm-main1}, but cannot be conjugate.

  Example~7.6 in \cite{DHL2016a} gives another variation on a product action, where the finite group $G$ is assumed to be simple, and one chooses a non-trivial proper subgroup $K \subset G$. The space $X \cong G/K$ is given   the left action by $G$, then using a product action  as above, we   obtain a minimal equicontinuous Cantor action  which is  locally quasi-analytic and has      discriminant  $\cD_x \cong K$.  
  If $K, K' \subset G$ are non-isomorphic proper subgroups with the same order, then  for product actions with the actions of $G$ and $G/K$ and $G/K'$ as above, we   obtain minimal equicontinuous Cantor actions which are  locally quasi-analytic  and continuously orbit equivalent,   but cannot be conjugate as their discriminants are not isomorphic.    
  
}
  \end{ex}
   
   \smallskip
    
  \begin{ex}[Semi-direct products]\label{ex-semidirect} 
 {\rm 
 
The direct product construction above can be made somewhat more interesting by using a semi-direct product construction, as in:

$\bullet$ ~ Example~7.5 in  \cite{DHL2016a}, where $G$ is the infinite dihedral group, and the discriminant $\cD_x \cong \mZ/2\mZ$. This action is non-homogeneous but is stable.

 $\bullet$ ~ Examples~8.8 and  8.9 in \cite{DHL2016a}, where $G$ is a generalized dihedral group, and the discriminant of the action is a Cantor subgroup.

 $\bullet$ ~ Example~7.1 in  \cite{DHL2016b}, where $G$ is a semi-direct product of $\mZ^n$ with a finite   subgroup $H  \subset \Sigma(n)$ of the permutation group on $n$ letters. 
 
Note that in all of these cited examples, the semi-direct product construction results in an extension of an infinite group by a finite normal subgroup.

}
  \end{ex}

\smallskip

    \begin{ex}[Full group modifications]\label{ex-modifications} 
    {\rm 

  Let  $(\fX,G,\Phi)$ be a minimal equicontinuous Cantor action. Assume the action $\Phi$ is topologically free, so that we can identify $G$ with its image $\Phi(G) \subset \Homeo(\fX)$. The idea of this construction is to modify the group $G$ by the addition of a finite set of  elements from its topological full group   $[[\fX,G,\Phi]]$. 
  
 Let $U \subset \fX$ be a proper adapted clopen subset, with stabilizer group $G_U \subset G$. Let $A \colon G_U \to G_U$ be a non-trivial automorphism. Define a new action $\Phi_U^A \colon G_U \to \Homeo(\fX)$ by setting, for $g \in G_U$:
 \begin{equation}\label{eq-twisted}
\Phi_U^A(g)(x) = A(g) \cdot x ~ {\rm for}~  x \in U, ~ {\rm and} ~  \Phi_U^A(g)(x) = x ~ {\rm for}~ x \in \fX - U \ .
\end{equation}

Let $H \subset \Homeo(\fX)$ denote the group generated by the images $\Phi(G)$ and $\Phi_U^A(G_U)$. 
Then $H$ is finitely generated, and 
by construction, the action $\Psi \colon H \times \fX \to \fX$ is continuously orbit equivalent to the given action $\Phi$.  Note that the action of $H$ is locally quasi-analytic  but not  topologically free.
  
 The clopen set $U$ is adapted to both actions, and the restriction of $\Psi$ to $U$ is just the conjugated action $\Phi_U^A \colon G_U \to \Homeo(\fX)$. Thus, the images $\Phi_U (G_U) \subset \Homeo(U)$ and $\Phi_U^A(G_U) \subset \Homeo(U)$ are equal. It follows that the actions are return equivalent. However,   $H$ is almost surely not isomorphic to $G$, so the actions $(\fX,G,\Phi)$ and $(\fX,H,\Psi)$ cannot be conjugate.

}
  \end{ex}

\smallskip
  
 \begin{ex}[Heisenberg]\label{ex-nilpotent}
{\rm 

In her thesis \cite{Dyer2015}, Dyer constructed examples of subgroup chains in the $2$-dimensional Heisenberg group, so that the discriminant of the Cantor action constructed from the chain is a Cantor group. These are probably the simplest examples of   Cantor actions 
by torsion-free finitely generated nilpotent groups. We recall  Example~8.5 from \cite{DHL2016a}.

Let $\cH$ be the discrete Heisenberg group, presented in the form $\cH = (\mZ^3, *)$ with the group operation $*$ given by $(x,y,z)*(x',y',z')=(x+x',y+y',z+z'+xy')$.
  Let $A_n =  \left( \begin{array}{cc} p^n & 0 \\ 0 & q^n\end{array}\right)$, where $p$ and $q$ are distinct primes, and consider the action represented by a group chain 
  $$G_0 =\cH ~ , ~ \{G_n \mid n \geq 1\} ~ {\rm where}~ G_n \equiv A_n \mZ^2 \times p^n \mZ  \ .$$
Then the discriminant   $\cD_x$ is a Cantor group, by explicit calculation \cite[Example~5.14]{Dyer2015}.
  It would be interesting to have a procedure for constructing nilpotent actions of higher dimension whose discriminants are Cantor groups.
 All such examples must be locally quasi-analytic by Corollary~\ref{cor-nilLQA}. 
}
  \end{ex}

\smallskip

    \begin{ex}[Actions on   rooted  trees]\label{ex-arboreal}
  {\rm 
Every minimal equicontinuous Cantor action can be interpreted as an action on a rooted spherically homogeneous tree. For actions as discussed in the previous examples, this fact is just one facet of their study, but for other actions it is a part of their definition. We briefly recall some basic facts about trees, then describe two classes of ``arboreal actions'' currently under investigation.

Let $T$ be a tree. That is, $T$ consists of the set of vertices $V = \cup_{\ell \geq 0} V_{\ell}$, where each $V_{\ell}$ is a finite set, and of the set of edges $E$, where an edge $t = [a,b]$ can join two vertices $a$ and $b$ only if $a \in V_{\ell}$ and $b \in V_{\ell+1}$, and every vertex $b \in V_{\ell +1}$ is joined to precisely one vertex in $V_{\ell}$. We assume that $V_0$ is a singleton, so $T$ is \emph{rooted} with root vertex $V_0 = \{v_0\}$. 
The  tree $T$ is \emph{spherically homogeneous}, if  there is a sequence of positive integers $(n_1, n_2, \ldots)$ such that for every $\ell \geq 1$, every   vertex $v \in V_{\ell-1}$ is joined by an edge to precisely $n_{\ell}$ vertices in $V_{\ell}$. We assume that $n_{\ell} > 1$ for all $\ell \geq 1$.

  An automorphism $g \in \Aut(T)$ of the rooted tree $T$ permutes the vertices within each level $V_{\ell}$, while preserving the connectedness of the tree. That is, for all $\ell \geq 1$ the vertices $v_{\ell-1} \in V_{\ell-1}$ and $v_{\ell} \in V_{\ell}$ are joined by an edge if and only if $g\cdot v_{\ell-1} \in V_{\ell-1}$ and $g \cdot v_{\ell} \in V_{\ell}$ are joined by an edge. Thus $\Aut(T)$ acts on infinite paths in the tree. 
  Here, a path in $T$ is an ordered infinite sequence of vertices $(v_\ell)_{\ell \geq 0}$ such that for all $\ell \geq 0$ the vertex $v_{\ell}$ is at level $V_{\ell}$, and the consecutive vertices $v_{\ell}$ and $v_{\ell+1}$ are joined by an edge.  Denote by $\cP$ the space of all infinite paths in $T$ with cylinder set topology. As $n_{\ell} > 1$ for   $\ell \geq 1$, the space of paths $\cP$ is a Cantor space. 

Let $(\fX, G, \Phi)$ be an effective minimal equicontinuous Cantor actions. The ``tree model'' for this action is constructed as follows. Choose a basepoint $x \in \fX$, and choose an adapted neighborhood system $\cU = \{U_{\ell} \mid \ell \geq 1\}$ at $x$, with stabilizer groups $G_{\ell} = G_{U_{\ell}}$. Set the root vertex $v_0 = \{\fX\}$, and let  the vertices at level $\ell \geq 1$ be the collection   of clopen sets $V_{\ell} \equiv \{g \cdot U_{\ell} \mid g \in G/G_{\ell}\}$. For $g,g' \in G$, there is an edge $t = [g \cdot U_{\ell}, g' \cdot U_{\ell +1}]$ between $g \cdot U_{\ell}$ and $g' \cdot U_{\ell +1}$ exactly when $g' \cdot U_{\ell +1} \subset g \cdot U_{\ell}$. 
This defines a spherically homogeneous tree $T_{\cU}$ where the  degree $n_{\ell} = [G_{\ell -1} : G_{\ell}]$.  The map  $\tau_{\infty} \colon  X_{\infty} \to \fX$ introduced in Section~\ref{sec-stable} can be interpreted as defining a homeomorphism 
$\tau_{\cU} \colon  \cP_{\cU} \to \fX$ which commutes with the action of $G$.
 
Now let $T$ be a rooted spherically homogeneous tree. 
It was shown in the work    \cite{BOERT1996}
 that $\Aut(T)$ is the infinite wreath product of the symmetric groups $\fS(n_{\ell})$, and so has a dense countably generated subgroup $G_0 \subset \Aut(T)$. This is  described in \cite[Section~2.5]{Lukina2018}, and see also Nekrashevych \cite{Nekrashevych2005}. Moreover, it is shown in \cite{Lukina2018} that the action of $G_0$ on the path space $\cP$ is not locally quasi-analytic. Here is a basic question:
\begin{prob}\label{prob-arboreal}
Let $T$ be  a rooted spherically homogeneous tree.  
Construct finitely-generated subgroups $G_0 \subset \Aut(T)$ such that the induced action on the  Cantor space of paths $\cP$ is not locally quasi-analytic.
\end{prob}
One interest in constructing finitely-generated groups acting on Cantor spaces is that each such action can be realized as the global monodromy action of a suspension foliation  as described in Section~8.1 of \cite{HL2017a}. If the action is a minimal equicontinuous Cantor action, then the resulting foliated space is homeomorphic to a weak solenoid as described in Section~\ref{sec-nil} above. 

The non-Hausdorff examples of Grigorchuk and Nekrashevych provide one class of  solutions to Problem~\ref{prob-arboreal}.
 The construction by Grigorchuk of his celebrated groups with intermediate growth in \cite{Grigorchuk1984} introduced the notion of groups acting on rooted spherically homogeneous trees which are generated by automata. The branch groups of \cite{BGS2003,Grigorchuk2000} are a far reaching extension of Grigorchuk's original construction.   It is easy to see that the action of the standard Grigorchuk group on the Cantor set of paths is not locally quasi-analytic. 
Nekrashevych studied in \cite{Nekrashevych2015,Nekrashevych2016} the dynamical properties of these actions, and showed that they often contain non-Hausdorff elements in their action groupoids. An action containing a non-Hausdorff element cannot be locally quasi-analytic by Proposition~\ref{prop-hausdorff2}.  

  } 
  \end{ex}
    \begin{ex}[Absolute Galois actions]\label{ex-galois}
  {\rm Odoni  initiated in \cite{Odoni1985} the study of arboreal representations of the absolute Galois group. Such representations are given by the action of a profinite group on a spherically homogeneous rooted tree. This field of study is   surveyed by   Jones  in \cite{Jones2013}. 
The second author    developed in \cite{Lukina2018}  a method of associating to an arboreal representation   an infinite chain of discrete groups, and gave examples of arboreal representations, some of  which are locally quasi-analytic, and some of which are not.
The examples of not locally quasi-analytic actions   in \cite{Lukina2018}  are by infinitely-generated groups, and it is an open problem to construct arboreal representations which are finitely generated and not locally quasi-analytic.

  } 
  \end{ex}

%%%%%%%%%%%%%%%%%%%%%%%%%%%%%%%%%%%%%%%%%%%%%%%%%%%%%%%

\end{document}